\documentclass{amsart}
\usepackage{enumerate}
\usepackage{amssymb}
\newtheorem{Theorem}{Theorem}[section]
\newtheorem{Definition}[Theorem]{Definition}
\newtheorem{Proposition}[Theorem]{Proposition}

\newtheorem{Lemma}[Theorem]{Lemma}
\newtheorem{Corollary}[Theorem]{Corollary}
\theoremstyle{remark}

\newtheorem{Example}[Theorem]{Example}

\def\il{\int_}
\def\hvi{\varphi}
\def\eps{\varepsilon}

\def\ovr{\overline}

\def\al{\alpha}

\def\th{\theta}

\def\Dl{\Delta}
\def\dl{\delta}

\def\bd{\partial}
\def\lm{\lambda}

\def\sm{\setminus}
\def\sbs{\subset}

\def\nea{\nearrow}
\def\sea{\searrow}

\def\supp{\operatorname{supp}}

\def\rk{\mathcal}
\def\re{{\mathbf {Re\,}}}
\def\im{\mathbf {Im\,}}
\def\be{\begin{enumerate}}
\def\ee{\end{enumerate}}
\def\bT{\begin{Theorem}}
\def\eT{\end{Theorem}}
\def\bP{\begin{Proposition}}
\def\eP{\end{Proposition}}
\def\bD{\begin{Definition}}
\def\eD{\end{Definition}}
\def\bE{\begin{Example}}
\def\eE{\end{Example}}
\def\bL{\begin{Lemma}}
\def\eL{\end{Lemma}}
\def\bC{\begin{Corollary}}
\def\eC{\end{Corollary}}

\def\oD{\ovr{\mathbb D}}
\def\aD{\mathbb D}

\def\aT{\mathbb T}
\def\aC{\mathbb C}

\def\E{{\mathcal E}}

\begin{document}
\title{Boundary values properties of functions in weighted Hardy spaces}
\author{Khim R. Shrestha}
\begin{abstract} In this paper we study the boundary values of harmonic and  holomorphic functions in the weighted Hardy spaces on the unit disk $\aD$. These spaces were introduced by Poletsky and Stessin in \cite{PS} for plurisubharmonic functions on hyperconvex domains $D\sbs\aC^n$ as generalizations of classical Hardy spaces. We show that in the case when $D$ is the unit disk $\aD$ the theory of boundary values for functions in these spaces is analogous to the classical one.
\end{abstract}
\address{Department of Mathematics,  Syracuse University, \newline
215 Carnegie Hall, Syracuse, NY 13244} \email{}
\maketitle

\section{Introduction}      %%%%%%%%%%%%%%%%%    SECTION 1      %%%%%%%%%%%%%%%%%%%%
\par In this paper we study the boundary values of harmonic and  holomorphic functions in the weighted Hardy spaces on the unit disk $\aD$. These spaces were introduced by Poletsky and Stessin in \cite{PS} for plurisubharmonic functions on hyperconvex domains $D\sbs\aC^n$ as generalizations of classical Hardy spaces. They are  parameterized by continuous negative plurisubharmonic exhaustion functions $u$ on $D$ and are denoted by $H^p_u(D)$. It was proved in \cite{PS} that $H^p_u(D)\sbs H^p(D)$ for all exhausting functions $u$.
\par As an example in Section \ref{Ex} shows that, in general, $H^p_u(\aD)\ne H^p(\aD)$. However, if $f\in H^p_u(\aD)$ then it belongs to $H^p(\aD)$ and, consequently, has radial boundary values $f^*$. The classical theory states that the Hardy norm of $f$ coincides with the norm of $f^*$ in $L^p(\lm)$, where $\lm$ is the normalized Lebesgue measure. Most of this paper is devoted to establishing analogous results for $H^p_u(\aD)$.
\par The definition of spaces $H^p_u(D)$ uses the measures $\{\mu_{u,r}\}$, $r<0$, (see Section \ref{basic}) introduced by Demailly in  \cite{D1}. These measures converge weak-$*$ in $C^*(\overline{\aD})$ to a positive measure $\mu_u$ supported by $\aT=\bd\aD$. As we show the measure $\mu_u$ replaces $\lm$ in the results about the spaces $H^p_u(\aD)$.
\par In Section \ref{bv} we define the Hardy spaces $h^p_u(\aD),\,p>1,$ of harmonic functions and prove that the norm of a function $h\in h^p_u(\aD)$ coincides with the norm of $h^*$ in $L^p(\mu_u)$. Also in this section we establish absolute continuity of $\mu_u$ with respect to $\lm$ and provide a formula for $\mu_u$.
\par In Section \ref{bvhar} for a function $h\in h^p_u(\aD)$ we show that the measures $h\mu_{u,r}$ converge weak-$*$ to the measure $h\mu_u$. This allows us to prove that $h$ has boundary values with respect to measure $\mu_u$ in the sense of \cite{Poletsky1}. After that in Section \ref{bvhol} we prove that the norm of a function $f\in H^p_u(\aD)$ coincides with the norm of $f^*$ in $L^p(\mu_u)$.
\par In section \ref{poh} we prove that the closed balls in $H^p_u(\aD)$ are closed in $H^p(\aD)$  and the space $H^p_u(\aD)$ is isometrically isomorphic to $H^p(\aD)$.
\par The author would like to thank Professor Evgeny Poletsky, who introduced this problem to him and constantly supported the development and improvement of this paper.

%%%%%%%%%%%%%%%%     SECTION 2       %%%%%%%%%%%%%%%%%%%%

\section{Basic facts} \label{basic}

\par Let $\aD$ be the unit disc $\{|z|< 1\}$ in $\aC$.  A continuous subharmonic function $u:\mathbb{D}\to [-\infty, 0)$ such that $u(z)\to 0$ as $|z|\to 1$ is called an exhaustion function. Following \cite{D1} for $ r < 0 $ we set
\[B_{u,r}  =\{z\in \mathbb{D}:u(z) < r\}\text{ and }
S_{u,r} = \{z\in \mathbb{D}: u(z) = r\}.\]
As in \cite{D1} we let $u_r = \max\{u,r\}$ and define the  measure
\[\mu_{u,r} = \Dl u_r -\chi_{\aD \setminus B_r}\Dl u, \] where $\Dl$ is the Laplace operator.
Clearly $\mu_{u,r}\ge 0$ and is  supported by $S_{u,r}$.
\par Let us denote by $\E$ the set of all continuous negative subharmonic exhaustion functions $u$ on $\aD$ such that
\[\il \aD\Dl u<\infty.\] In the same paper Demailly (see Theorems 1.7 and 3.1 there) proved the following result which we adapt to the case of $\aD$.
\bT [Lelong--Jensen formula]\label{jlf}
Let $\phi$ be a subharmonic function on $\mathbb{D}$. Then $\phi$ is $\mu_{u,r}$-integrable for every $r < 0$ and \[\mu_{u,r}(\phi) = \int_{B_{u,r}} \phi\, \Dl u + \int_{B_{u,r}}(r-u)\,\Dl\phi.\]
Moreover, if $u\in\E$ then the measures $\mu_{u,r}$ converge weak-$*$ in $C^*(\oD)$  to a measure $\mu_u\ge 0$ supported by $\aT$ as $r\to 0^-$.
\eT
He also derived from this theorem the following
\bC\label{C:imu} If $\phi$ is a non-negative subharmonic function, then the function $r\to \mu_{u,r}(\phi)$ is increasing on $(-\infty, 0)$.\eC
\par Using the measures $\mu_{u,r}$ Poletsky and Stessin introduced the Hardy spaces associated with an exhaustion $u\in\E$. For $0 < p < \infty$ we define the space $H^p_u(\mathbb{D})$ consisting of the functions $f(z)$ analytic in $\mathbb{D}$ and satisfying \[\|f\|_{H^p_u}^p = \varlimsup_{r\to 0^-}\il {S_{u,r}} |f|^p\, d\mu_{u,r} < \infty.\] By Corollary \ref{C:imu} we can replace the $\varlimsup$ in the above definition with $\lim$.  By Theorem \ref{jlf} and the monotone convergence theorem it follows that,
\begin{equation}\label{nf} \|f\|^p_{H^p_u} = \il {\aD} |f|^p\,\Dl u - \il {\aD} u\,\Dl |f|^p.
\end{equation}
The classical Hardy spaces correspond to $u(z)=\log|z|$ (see Section 4 in \cite{PS}) and will be denoted by $H^p(\aD)$.

\par It was proved in \cite{PS} that:\be
\item the spaces $H^p_u(\aD)$ are Banach when $p\ge 1$ (Theorem 4.1);
\item if $v,u\in\E$ and $v\le u$ on $\aD$, then $H^p_v(\aD)\sbs H^p_u(\aD)$ and if $f\in H^p_v(\aD)$ then $\|f\|_{H^p_u}^p\le \|f\|_{H^p_v}^p$.
\ee
\par Thus by Hopf's lemma the space $H^p_u(\mathbb{D})$ is contained in the classical Hardy space $H^p(\mathbb{D})$.
\section{Example}\label{Ex}    %%%%%%%%%%%%%%%%   SECTION 3   %%%%%%%%%%%%%%%%%%

\par Having known that the space $H^p_u(\aD)$ is contained in $H^p(\aD)$, a question arises naturally whether $H^p_u(\aD)$ is properly contained in $H^p(\aD)$ or $H^p_u(\aD)$ can also be equal to $H^p(\aD)$. In \cite{PS} it has been proved that if $\Dl u$ is compactly supported then $H^p_u(\aD) = H^p(\aD)$. However, this is not the case in general. Now we construct a subharmonic function $u(z)\in\E$ on $\aD$ for which $H^2_u(\aD)\ne H^2(\aD)$.

\bL\label{construction} If $0 <\beta < 1$ the integral
\[\int_0^1 \log\left|\frac{s-t}{1 - ts}\right|\frac{ds}{(1-s)^\beta},\,\,\,0 < t < 1,\]
tends to $0$ as $t\to 1$. \eL

\begin{proof}
Write
\begin{align*}
\int_0^1 \log\left|\frac{s-t}{1 - ts}\right|\frac{ds}{(1-s)^\beta} & = \int_0^t \log\left(\frac{t-s}{1 - ts}\right)\frac{ds}{(1-s)^\beta} + \int_t^1 \log\left(\frac{s-t}{1 - ts}\right)\frac{ds}{(1-s)^\beta} \\
& =\text{ I }+ \text{ II }.
\end{align*}
Make a substitution of $\displaystyle s=\frac{x+t}{1+tx}$
in II to get
\begin{align*}
\text{ II } & =  (1+t)(1-t)^{1-\beta}\int_0^1 \frac{\log x}{(1-x)^\beta (1 + tx)^{2-\beta}}\,dx\\
& \ge (1+t)(1-t)^{1-\beta}\int_0^1 \frac{\log x}{(1-x)^\beta}\,dx\\
& \to 0\,\,\, \text{ as $ t\to 1$ when $ 0 < \beta < 1$. }
\end{align*}
Again, make substitution of $\displaystyle  s=\frac{t-x}{1 - tx}$ in I to get
\begin{align*}
\text{ I } & = (1+t)(1-t)^{1-\beta}\int_0^t \frac{\log x}{(1+x)^\beta (1-tx)^{2-\beta}}\,dx\\
& \ge (1+t)(1-t)^{1-\beta}\int_0^t \frac{\log x}{(1-tx)^{2-\beta}}\,dx\\
& = t(1+t)(1-t)^{1-\beta}\int_0^1 \frac{\log(tx)}{(1-t^2x)^{2-\beta}}\,dx\\
& \ge t(1+t)(1-t)^{1-\beta} \left( \int_0^1\frac{\log t}{(1-t^2x)^{2-\beta}}\,dx + \int_0^1\frac{\log x}{(1-x)^{2-\beta}}\,dx\right)\\
& \to 0 \text{ as $t \to 1$ when $0 < \beta < 1$.}
\end{align*}
Thus $u(t) \to 0$ as $t \to 1$ when $0 < \beta < 1$.
\end{proof}

\par Now define a function $u(z): \aD \to [-\infty,0)$ by
\[ u(z) = \int _0^1\log\left|\frac{z-s}{1-sz}\right|\frac{ds}{(1-s)^\beta},\]
where $\beta$ is a number between $0$ and $1$.
The function $u(z)$ is subharmonic. If $z, w \in \aD$, then by the inequality (see \cite[Lemma 4.5.7]{Ransford})
\[\left|\frac{|z|-|w|}{1-|w||z|}\right| \le \left|\frac{z-w}{1-\bar{w}z}\right|\]
and Lemma \ref{construction} it follows that $u(z)\to 0$ as $|z| \to 1$. Also \[ \il \aD\Dl u = \int_0^1 \frac{dx}{(1-x)^\beta}< \infty.\]
 Thus $u\in \E$ and, for this $u$, we show that $H^2_u(\aD) \ne H^2(\aD)$.

\bT For $\frac{1-\beta}{2} \le \al< \frac{1}{2}$ the function \[f(z) = \frac{1}{(1-z)^\al}\] is in $H^2(\aD)$ but not in $H^2_u(\aD)$.\eT
\begin{proof}
The function $f(z) = \frac{1}{(1-z)^\al}$ belongs to $H^p(\aD)$ for every $\al < \frac{1}{p}$  (see \cite{Duren1}, page 78). Hence $f(z)\in H^2(\aD)$ for $\al < \frac{1}{2}$. On the other hand, by \eqref{nf}
\[ \|f\|^2_{H^2_u} \ge \il \aD|f|^2\,\Dl u = \int_0^1 \frac{1}{(1-x)^{2\al +\beta}}\,dx = \infty\]
when $2\al +\beta \ge 1$.  Hence $f(z) \notin H^2_u(\aD)$ for $\al \ge \frac{1-\beta}{2}$.

\end{proof}

%%%%%%%%%%%%%%    SECTION 4   %%%%%%%%%%%%%%%%
\section{The Hardy spaces of harmonic functions and the measure $\mu_u$}\label{bv}

\par Let us denote by $h^p_u(\aD),\, p > 1$, $u\in\E$, the space of harmonic functions $h $ on $\aD$ such that
\[ \|h\|^p_{u,p} = \lim_{r\to 0^-} \int_{S_{u,r}} |h|^p\, d\mu_{u,r} < \infty.\]
By Corollary 3.2 in \cite{PS}, $h_u^p(\aD)\sbs h^p(\aD)$. Thus if $h\in h^p_u(\aD)$, then $h$ has radial boundary values $h^*$ on $\bd\aD$. We have the following theorem.
\bT \label{hbvalue}Let $h\in h^p_u(\aD),\, p > 1$. Then
\[ \|h\|^p_{u,p} = \int_{\aT} |h^*(e^{i\theta})|^p\, d\mu_u(\th).\]\eT

\begin{proof} Let $\lm$ be the normalized Lebesgue measure on $\aT$. The least harmonic majorant on $\aD$ of the subharmonic function $|h|^p$ is the Poisson integral of $|h^*|^p$. By the Riesz Decomposition Theorem
\[|h(w)|^p=\il {\aT}|h^*(e^{i\th})|^pP(w,e^{i\th})\,d\lm(\th) +
\il {\aD}G(w,z)\Dl|h|^p(z),\]  where $P$ is the Poisson kernel and $G$ is the Green kernel.
\par By Lelong--Jensen formula and the monotone convergence theorem we have
\[\|h\|^p_{u,p}=\il {\aD}|h|^p\Dl u-\il {\aD}u\Dl|h|^p.\]
Again by the Riesz formula,
\begin{equation}\label{rf} u(z)=\il {\aD}G(z,w)\Dl u(w). \end{equation}
Hence, by Fubini--Tonnelli's  Theorem and the symmetry of the Green kernel
\[\il {\aD}u(z)\Dl|h|^p(z)=
\il {\aD}\left(\il {\aD}G(w,z)\Dl|h|^p(z)\right)\Dl u(w)\]
and
\begin{equation}\begin{aligned}
\|h\|^p_{u,p} &=\il {\aD}\left(|h(w)|^p-\il {\aD}G(w,z)\Dl|h|^p(z)\right)\Dl u(w)\notag\\
&=\il {\aD}\left(\il{\aT}|h^*(e^{i\th})|^pP(w,e^{i\th})\,d\lm(\th)\right)\Dl u(w)\\
&=\il{\aT}\left(\il{\aD}P(w,e^{i\th})\Dl u(w)\right)|h^*(e^{i\th})|^p\,d\lm (\th).\notag
\end{aligned}\end{equation}
\par Let
\begin{equation}\label{alpha}\al(e^{i\th})=\il{\aD}P(w,e^{i\th})\Dl u(w).\end{equation}
Then
\[\|h\|^p_{u,p}=\il{\aT}|h^*(e^{i\th})|^p\al(e^{i\th})\,d\lm(\th).\]
\par Let $\phi$ be a continuous function on $\aT$ and let $h$ be its harmonic extension to $\aD$. Then $h^*=\phi$ and by Theorem \ref{jlf}
\[\|h\|^p_{u,p}=\il{\aT}|\phi(e^{i\th})|^p\,d\mu_u(\th).\]
Hence $\mu_u=\al\lm$ and $\al\in L^1(\lm)$.
Consequently, for any $h\in h^p_u(\aD)$
\[\|h\|^p_{u,p}=\il{\aT}|h^*(e^{i\th})|^p\,d\mu_u(\th).\]
\end{proof}

\par We collect the information about the measure $\mu_u$ in the following proposition.

\bP\label{abscont} The measure $\mu_u=\al\lm$, where the function $\al(e^{i\th})$ has the following properties:
\begin{enumerate}[(i)]
\item $\al(e^{i\th})\in L^1(\lm)$.
\item $\al(e^{i\th}) = \il {\aD} P(z, e^{i\th})\,\Dl u(z)$.
\item $\al(e^{i\th})$ is lower semicontinuous.
\item $\al(e^{i\th}) \ge c > 0$ on $\aT$.
\item $\al(e^{i\th})$ need not to be necessarily bounded. 
\end{enumerate}
\eP
\begin{proof} Everything except $(iii),\,(iv)\, \text{ and } (v)$ follow from the proof of the above theorem. Let $e^{i\th_j}\to e^{i\th_0}$ in $\aT$. By Fatou's lemma
\[ \liminf_{j\to\infty} \al(e^{i\th_j}) = \liminf_{j\to\infty} \il {\aD} P\left(z,e^{i\th_j}\right)\Dl u(z) \ge \il {\aD} P\left(z,e^{i\th_0}\right) \Dl u(z) = \al(e^{i\th_0}).\] This proves $(iii)$. 
\par Let $v(z) = \log|z|$. By Hopf's lemma there is a constant $c>0$ such that $cu(z)<v(z)$ near $\aT$.
It follows from \cite[Theorem 3.8]{D1} that $\mu_v\le c\mu_u$. Since $\mu_v = \lm$, $(iv)$ follows.
\par
For the exhaustion function  constructed in Section \ref{Ex},
\[ \int_\aD P(z, 1) \Dl u = \int_0^1 \frac{1+x}{1-x}\cdot \frac{1}{(1-x)^\beta}\, dx = \infty\]
when $\beta > 0$. This proves $(v)$. 
\end{proof}

\par In the proof of the theorem \ref{hbvalue} we have deduced the norm of the functions $h\in h^p_u(\aD),\, p>1$ to 
\[\|h\|^p_{u,p} = \il {\bd\aD}\left(\il \aD P(w, e^{i\th})\,\Dl u(w)\right)|h^*(e^{i\th})|^p\,d\lm.\]
Since $\frac{\bd }{\bd n}G(z, w)|_{z=e^{i\th}}=P(e^{i\th}, w)$, from the Riesz formula (\ref{rf})  we get
\[ \frac{\bd u}{\bd n}(e^{i\th}) = \il \aD P(w, e^{i\th})\,\Dl u(w)\]
and therefore the norm can be written as 
\[\|h\|^p_{u,p} = \il {\bd\aD}\frac{\bd u}{\bd n}(e^{i\th})|h^*(e^{i\th})|^p\,d\lm.\]

From this deduction it is clear that if $u\in\E$ is such that $\frac{\bd u}{\bd n}(e^{i\th})$ is bounded then $h^p_u(\aD) = h^p(\aD),\,p>1$.

%%%%%%%%%%%%%%% Section 5 %%%%%%%%%%%%%%%%%%%%

\section{Boundary values of harmonic functions with respect to\\ the measures $\mu_{u,r}$}\label{bvhar}
\par While functions in $h^p_u(\aD)$, $p>1$, have radial limits $\mu_u$-a.e., we are interested in the analogs of more subtle classical properties of boundary values. For example, if $h\in h^p(\aD)$ then it is known that the measures $h(re^{i\th})\lm(\th)$ converge weak-$*$ in $C^*(\aT)$ to $h(e^{i\th})\lm(\th)$ as $r\to1^-$.
\par In this section we will establish the analogs of these statements.
\bT \label{wlimit} Let $h \in h^p_u(\aD),\, p > 1$. Then the measures $\{h\mu_{u,r}\}$ converge weak-${\ast}$ to $h^*\mu_u$ in $C^*(\oD)$ when $r\to0^-$. \eT
\begin{proof}
Since the space $C(\oD)$ is separable the weak-$*$ topology on the balls in $C^*(\oD)$ is metrizable. Thus it suffices to show that for any sequence $r_j\nea0$ and any $\phi\in C(\oD)$ we have
\[\lim_{j\to\infty}\il{S_{u,r_j}}\phi h\,d\mu_{u,r_j}=\il {\bd\aD}\phi h^*\,d\mu_u.\]
\par We introduce functions
\[p_r(e^{i\th})=\il{S_{u,r}}P(z,e^{i\th})\,d\mu_{u,r}(z)=\il{B_{u,r}}P(z,e^{i\th})\,\Dl u(z),\] where the last equality follows from Theorem \ref{jlf} because $\Dl h\equiv0$. Hence $p_r(e^{i\th})\nea\al(e^{i\th})$.
\par  Due to the uniform continuity of $\phi$ and the formula for $P(z,e^{i\th})$,  for every $\th\in[0,2\pi]$ and for every $\eps>0$ there is $\dl>0$ such that $|P(z,e^{i\th})|<\eps$ when $z$ is close to boundary and $|z-e^{i\th}|>\dl$ and $|\phi(z)-\phi(e^{i\th})|<\eps$ when $|z-e^{i\th}|\le\dl$. Hence, when $r$ is sufficiently close to 0,
\begin{equation}\begin{aligned}
&\left|\il { S_{u,r}}\phi(z)P(z,e^{i\th})\,d\mu_{u,r}(z)-
\il{S_{u,r}}\phi(e^{i\th})P(z,e^{i\th})\,d\mu_{u,r}(z)\right| \notag\\
\le&\il{S_{u,r}\sm\oD(e^{i\th},\dl)}|\phi(z)-\phi(e^{i\th})|P(z,e^{i\th})\,d\mu_{u,r}(z)\notag\\
& +\il{S_{u,r}\cap\oD(e^{i\th},\dl)}|\phi(z)-\phi(e^{i\th})|P(z,e^{i\th})\,d\mu_{u,r}(z) \notag\\
\le&2M\eps+\eps p_r(e^{i\th}),\notag
\end{aligned}\end{equation}
where $M$ is the uniform norm of $\phi$ on $\ovr{\aD}$.
\par Now,
\begin{equation}\begin{aligned}\il{S_{u,r}}\phi(z)h(z)\,d\mu_{u,r}(z)=&
\il{S_{u,r}}\phi(z)\left(\il{\aT}h^*(e^{i\th})P(z,e^{i\th})\,d\lm(\th)\right)\,d\mu_{u,r}(z)\notag\\=
&\il{\aT}h^*(e^{i\th})\left(\il{S_{u,r}}\phi(z)P(z,e^{i\th})\,d\mu_{u,r}(z)\right)\,d\lm(\th).\notag
\end{aligned}\end{equation}
Hence,
\begin{equation}\begin{aligned}
&\left|\il{S_{u,r}}\phi(z)h(z)\,d\mu_{u,r}(z)-
\il{\aT}\phi(e^{i\th})h^*(e^{i\th})\,d\mu_u(\th)\right|\notag\\
\le& \left|\il{S_{u,r}}\phi(z)h(z)\,d\mu_{u,r}(z)- \il {\aT} \phi(e^{i\th})h^*(e^{i\th})p_r(e^{i\th})\,d\lm(\th)\right|\notag\\
&+\left|\il {\aT} \phi(e^{i\th})h^*(e^{i\th})p_r(e^{i\th})\,d\lm(\th)-
\il{\aT}\phi(e^{i\th})h^*(e^{i\th})\,d\mu_u(\th)\right|\notag\\
=&\left|\il{\aT} h^*(e^{i\th})\left(\il{S_{u,r}}(\phi(z)-
\phi(e^{i\th}))P(z,e^{i\th})\,d\mu_{u,r}(z)\right)\,d\lm(\th)\right|\notag\\
&+\left|\il {\aT} \phi(e^{i\th})h^*(e^{i\th})\left(p_r(e^{i\th})-\al(e^{i\th})\right)\,d\lm(\th)\right|\notag\\
\le&\eps\il{\aT}\left |h^*(e^{i\th})\right|(2M+p_r(e^{i\th}))\,d\lm(\th)  + M\il {\aT} \left|h^*(e^{i\th})\right|\left|p_r(e^{i\th})-\al(e^{i\th})\right|\,d\lm(\th).\notag
\end{aligned}\end{equation}
 Now,
\begin{equation}\begin{aligned}
 \il{\aT}\left |h^*(e^{i\th})\right|(2M+p_r(e^{i\th}))\,d\lm(\th) \le &\il{\aT}\left |h^*(e^{i\th})\right|(2M+\al(e^{i\th}))\,d\lm(\th)\notag\\
\le & 2M\|h^*\|_{L^p} + \|h\|_{u,p}.\notag
\end{aligned}\end{equation}
Since $\left|p_r(e^{i\th})-\al(e^{i\th})\right| \sea 0$ and $\left|p_r(e^{i\th})-\al(e^{i\th})\right| < \al(e^{i\th})$ with $\left|h^*(e^{i\th})\right|\al(e^{i\th}) \in L^1(\lm)$, by the monotone convergence theorem,
\[\il {\aT} \left|h^*(e^{i\th})\right|\left|p_r(e^{i\th})-\al(e^{i\th})\right|\,d\lm(\th) \to 0\]
Thus, since $\eps$ is arbitraty,
\[\left|\il{S_{u,r}}\phi(z)h(z)\,d\mu_{u,r}(z)-
\il{\aT}\phi(e^{i\th})h^*(e^{i\th})\,d\mu_u(\th)\right|\to 0.\]
The proof is complete.
\end{proof}
\par In \cite{Poletsky1} Poletsky introduced the weak and strong limit values for a sequence $\{\phi_j\}$ of Borel functions defined on compact subsets $K_j$ of a compact set $K$ with respect to a sequence of regular Borel measures $\{\mu_j\}$ supported by $K_j$ and  converging weak-$\ast$ in $C^*(K)$ to a finite measure $\mu$. If the measures $\{\phi_j\mu_j\}$ converge weak-$\ast$ in $C^*(K)$ to a measure $\phi_\ast\mu$ , then the function $\phi_\ast$ is called the \emph {weak limit values} of $\{\phi_j\}$.
\par We say that the sequence $\{\phi_j\}$ has a \emph{strong limit values} on $\supp \mu$ with respect to $\{\mu_j\}$ if there is a $\mu$-measurable function $\phi^*$ on $\supp \mu$ such that for any $b > a$ and any $\epsilon, \dl > 0$ there is $j_0$ and an open set $O \subset K$ containing $G(a,b) = \{x\in \supp\mu: a\le \phi^*(x) < b\} $ such that
\[\mu_j(\{\phi_j < a-\epsilon\}\cap O) +\mu_j(\{\phi_j > b+ \epsilon\}\cap O) < \dl\]
when $j\ge j_0$. The function $\phi^*$ is called the \emph{strong limit values} of $\{\phi_j\}$.
\par Following the definition in \cite{Poletsky1}, we say that a function $h\in h^p_u(\aD)$ has \emph{ boundary values } with respect to the measures $\mu_{u,r}$ if it has strong limit values with respect to $\{\mu_{u,r_j}\}$ for any sequence $r_j\nea 0$ and these strong limit values do not depend on the choice of a sequence.

\bT \label{hslim} Let $h\in h^p_u(\aD),\, p>1$. Then $h$ has the boundary values equal to $h^*$ with respect to $\{\mu_{u,r}\}$. \eT

\begin{proof} Let $r_j$ be any increasing sequence of numbers converging to 0. By Theorem \ref{wlimit} the measures $h\mu_{u,r}$ converge weak-$*$ in $C^*(\oD)$ to the measure $h^*\mu_u$. By Theorem \ref{hbvalue}
\[\lim_{j\to\infty}\il{S_{u,r_j}}|h|^p\,d\mu_{u,r_j}=\il{\aT}|h^*|^p\,d\mu_u.\]
By \cite[Theorem 3.6]{Poletsky1} the sequence of the function $h|_{S_{u,r_j}}$ has the strong boundary values equal to $h^*$.
\end{proof}

%%%%%%%%%%%%%%%%%% Section 6 %%%%%%%%%%%%%%%%%%%%%%%%%

\section{Boundary values of analytic functions with respect to\\ the measures $\mu_{u,r}$}\label{bvhol}

\par In this section we prove results analogous to those in two previous sections but for $p>0$. To consider the Hardy spaces for $0<p\le1$ we need a factorization theorem.
\par From the classical theory we know that every function $f\in H^p(\aD),\, p>0,\, f\not\equiv 0$ can be factorized into $f(z) = \beta(z)g(z)$ where $\beta(z)$ is a Blaschke product with same zeros as $f$ and $g$ is a non-vanishing function in $H^p(\aD)$ with $\|g\|_{H^p} = \|f\|_{H^p}$. Let us show that the similar result holds for the functions in $H^p_u(\aD)$.
\bT\label{Bproduct}
Let $f(z)\in H^p_u(\aD),\,\, p > 0$ and $f(z)\not\equiv 0$. Then there exists a function $g(z) \in H^p_u(\aD)$, $g(z)\ne 0$ in $\aD$, such that
\[
f(z) = \beta(z)g(z)\,\, \text{ and }\,\,\|g\|_{H^p_u}=\|f\|_{H^p_u},
\]
where $\beta(z)$ is a Blaschke product having the same zeros as $f$.
\eT
\begin{proof} We mimic the proof of the classical version \cite[Theorem 2.3]{Garnett}.
Let $\{a_j\}$ be the zeros of $f(z)$ in $\aD$ not necessarily all distinct. We may assume that $a_j\ne 0$ for all $j$ since otherwise if $0$ is the zero of order $m$ then we write $f(z) = z^m\tilde{f}(z)$ and work with $\tilde{f}(z)$. Then
\[
\beta(z) = \prod_{j=1}^\infty\frac{-\overline{a}_j}{|a_j|}\frac{z-a_j}{1-\overline{a}_jz}.
\]
From classical theory  we have $g(z) = \frac{f(z)}{\beta(z)} \in H^p(\aD)$. We show that $g(z) \in H^p_u(\aD)$.

\par Write
\[
 g_N(z) = \frac{f(z)}{\beta_N(z)},\,\,\text{ where } \beta_N(z) = \prod_{j=1}^N\frac{-\overline{a}_j}{|a_j|}\frac{z-a_j}{1-\overline{a}_jz}.
\]
For fixed $N$, $|\beta_N(z)| \to 1$ uniformly as $|z| \to 1$. So for given $\eps > 0$ there exists $\rho_0 > 0$ such that $|\beta_N(z)| > 1 -\eps$ when $|z| > \rho_0.$
Thus near $\aT$ we have
\[
|g_N(z)| < \frac{|f(z)|}{1-\eps}.
\]
Since $\eps$ is arbitrary and $\mu_{u,r}(|f|^p)$ is an increasing function of $r$, it follows that
\[
\il {S_{u,r}}|g_N(z)|^p\,d\mu_{u,r} \le \|f\|^p_{H^p_u}.
\]
Since $|g_N(z)|\nea |g(z)|$, by the monotone convergence theorem,
\[
\il {S_{u,r}}|g(z)|^p\,d\mu_{u,r} =\lim_{N\to\infty}\il {S_{u,r}}|g_N(z)|^p\,d\mu_{u,r} \le\|f\|^p_{H^p_u}.
\]
Hence $\|g\|_{H^p_u} \le \|f\|_{H^p_u}$. The reverse inequality is trivial because $|f(z)|\le |g(z)|$ in $\aD$. Thus $\|g\|_{H^p_u} = \|f\|_{H^p_u}.$ This completes the proof.
\end{proof}
\par Since $H^p_u(\aD) \sbs H^p(\aD)$, any $f\in H^p_u(\aD)$ has radial limits $f^*(e^{i\th})\,\lm$-a.e. But it is not clear that $\|f\|_{H^p_u}\ge \|f^*\|_{L^p(\mu_u)}$. The theory of weak and strong limit values in \cite{Poletsky1} provides sufficient conditions for this estimate. To implement these conditions we have to show that the existence of strong limit values for $f\in H^p_u(\aD)$.
\bT\label{wlimit1}
Any function $f\in H^p_u(\aD),\, p>1,$ has the weak limit values equal to $f^*$ with respect to the measures $\{\mu_{u,r}\}$.
\eT
\begin{proof} Follows directly from Theorem \ref{wlimit}. \end{proof}
\bT \label{slim} Let $f\in H^p_u(\aD),\, p>1$. Then $|f|$ has the boundary values equal to      $|f^*|$ with respect to $\{\mu_{u,r}\}$. \eT
\begin{proof}
For $f\in H^p_u(\aD),\, \re f \text{ and } \im f\in h^p_u(\aD)$. Hence the corollary follows from Theorem \ref{hslim} and \cite[Theorem 3.3]{Poletsky1} by writing $|f|^2 = (\re f)^2 + (\im f)^2$.
\end{proof}

\par Now we prove the most important theorem of the section:
\bT\label{bvalue}
Let $f\in H^p(\aD),\,p >0 $. Then $f\in H^p_u(\aD)$ if and only if $f^*(e^{i\th})\in L^p(\mu_u)$. Moreover, $\|f\|_{H^p_u} = \|f^*\|_{L^p(\mu_u)}.$
\eT

\begin{proof}
First, we prove the theorem for $p > 1$. Let $f^*\in L^p(\mu_u)$. There exists $f_j^* \in C(\aT)$ such that
\[\|f^*_j-f^*\|_{L^p(\mu_u)} \to 0 \,\,\,\text{ as } j\to \infty.\]
By Proposition \ref{abscont},
\[\|f^*_j-f^*\|_{L^p(\lm)} \to 0 \,\,\,\text{ as } j\to \infty.\]
We know that $f(z)$ is the Poisson integral of its boundary value $f^*(e^{i\th})$ \cite[Theorem 3.1]{Duren}, that is,
\[f(z) = \int_0^{2\pi}P(z, e^{i\th}) f^*(e^{i\th})\,d\lm(\th). \]
If we take
\[f_j(z) = \int_0^{2\pi}P(z,e^{i\th}) f_j^*(e^{i\th})\,d\lm(\th)\]
by H\"{o}lder's inequality,
\begin{align*}
|f_j(z) -f(z)| & =  \left|\int_0^{2\pi} \left(f_j^*(e^{i\th}) -f^*(e^{i\th})\right)P(z,e^{i\th})\,d\lm(\th)\right|\\
 &   \le \left(\int_0^{2\pi}\left|f_j^*(e^{i\th}) -f^*(e^{i\th})\right|^p\,d\lm(\th)\right)^\frac{1}{p}\left(\int_0^{2\pi}P^q(z, e^{i\th})\,d\lm(\th)\right)^\frac{1}{q}.
\end{align*}
The last integral is, evidently, bounded on compact sets in $\aD$ and hence $f_j\to f$ uniformly on compacta. Therefore
\[
\lim_{j\to\infty} \int_{S_{u,r}}|f_j|^p\,d\mu_{u,r} = \int_{S_{u,r}}|f|^p\,d\mu_{u,r}.
\]
The weak-$*$ convergence of $\mu_{u,r}$ gives
\[
\lim_{r\to 0^-}\int_{S_{u,r}}|f_j|^p\,d\mu_{u,r} = \int_{\aT}|f_j|^p\,d\mu_u.
\]
Since  $f_j(z)$ is harmonic, $|f_j|^p$ is subharmonic and by Corollary \ref{C:imu},\, $\mu_{u,r}(|f_j|^p)$ is an increasing function of $r$. It follows, for each $j$, that
\[
\int_{S_{u,r}}|f_j|^p\,d\mu_{u,r}\le \int_{\aT}|f_j|^p\,d\mu_u =\int_{\aT}|f_j^*|^p\,d\mu_u .
\]
Hence
\[
\int_{S_{u,r}}|f|^p\,d\mu_{u,r} =\lim_{j\to \infty}\int_{S_{u,r}}|f_j|^p\,d\mu_{u,r}\le \lim_{j\to\infty}\int_{\aT}|f_j^*|^p\,d\mu_u = \int_{\aT}|f^*|^p\,d\mu_u.
\]
Therefore  $\|f\|_{H^p_u} \le \|f^*\|_{L^p(\mu_u)}$ and $f\in H^p_u(\aD)$.

\par Let $f\in H^p_u(\aD)$. Then by Corollary \ref{slim},  $|f|$ has the boundary values $|f^*|$ with respect to $\{\mu_{u,r}\}$. By \cite[Theorem 3.5]{Poletsky1}, it follows that
\[\|f^*\|_{L^p(\mu_u)}\le \|f\|_{H^p_u}.\]
Hence $f^* \in L^p(\mu_u)$ and $\|f\|_{H^p_u} = \|f^*\|_{L^P(\mu_u)}$.

\par Now we prove the theorem for $0 < p \le 1$.
Let $f\in H^p(\aD)$. Then we have the factorization $f(z)= \beta(z) g(z)$ where $\beta(z)$ is a Blaschke product and $g(z)$ is a non-vanishing function in $H^p(\aD)$. Suppose $f^*\in L^p({\mu_u})$. Since $|f^*| = |g^*|\,\lm$-a.e. (and hence  $\mu_u$-a.e.), $g^* \in L^p(\mu_u)$. It follows from the proof for $p>1$ and the fact that $g^\frac{p}{2} \in H^2(\aD)$ and $(g^*)^\frac{p}{2} \in L^2(\mu_u)$ that
\[
\|g^\frac{p}{2}\|_{H^2_u} \le \|(g^*)^\frac{p}{2}\|_{L^2(\mu_u)}.
\]
This implies
\[
\|g\|_{H^p_u} \le \|g^*\|_{L^p(\mu_u)}.
\]
Since $|f(z)|\le|g(z)|$ in $\aD$ we get
\[
\|f\|_{H^p_u} \le \|f^*\|_{L^p(\mu_u)}
\]
and hence $f\in H^p_u(\aD)$.

\par On the other hand if $f\in H^p_u(\aD)$ then by Theorem \ref{Bproduct},\, $f(z) = \beta(z)g(z)$ where $g(z)$ is non-vanishing function in $H^p_u(\aD)$. Since $g^{\frac{p}{2}}\in H^2_u(\aD),\, |g^{\frac{p}{2}}|$ has the boundary values $|(g^{\frac{p}{2}})^*|$ with respect to $\{\mu_{u,r}\}$. Then by  \cite[Theorem 3.5]{Poletsky1},
\[
\|(g^\frac{p}{2})^*\|_{L^2(\mu_u)} \le \|g^\frac{p}{2}\|_{H^2_u}.
\]
This implies
\[
\|g^*\|_{L^p(\mu_u)} \le \|g\|_{H^p_u}
\]
and hence
\[
\|f^*\|_{L^p(\mu_u)} \le \|f\|_{H^p_u}.
\]
Thus $f^*\in L^p(\mu_u)$ and $\|f\|_{H^p_u} = \|f^*\|_{L^p(\mu_u)}$.
\end{proof}

 %%%%%%%%%     SECTION 7    %%%%%%%%%%

\section{ Properties of $H^p_u(\aD)$}\label{poh}

Note that $H^p_u(\aD)$ is not a closed subspace of $H^p(\aD)$ because both spaces contain $H^\infty(\aD)$.  However, the closed balls in $H^p_u(\aD)$ are closed in $H^p(\aD)$.

\bT\label{unitball}
The closed unit ball
\[ B_{u,p}(1)= \left\{ f\in H^p_u(\aD): \|f\|_{H^p_u}\le 1\right\}\]
 in $H^p_u(\aD),\,p > 0$, is closed in $H^p(\aD)$.\eT

\begin{proof}
The case $p=\infty$ is obvious. Let $\left\{f_j\right\} \subset B_{u,p}(1)$ be such that $f_j \to f$ in $H^p(\aD)$, i.e.
\[\sup_{0\le r < 1}\int_0^{2\pi}\left|f_j(r e^{i\theta}) - f(r e^{i\theta})\right|^p\, d\lm(\th) \to 0\,\,\,\text{ as } j \to \infty.\]
By formula (3.2) in \cite{PS} if $|z|<r$ then
\[|f(z)-f_j(z)|^p\le\il{|w|=r}|f(re^{i\th})-f_j(re^{i\th})|^p\,d\lm(\th)\le\|f_j-f\|_{H^p}.\]
Hence the functions $f_j \to f $ uniformly on compacta.
\par Now
$$\int_{S_{u,r}}\left|f_j(z)\right|^p\,d\mu_{u, r} \to \int_{S_{u,r}}\left|f(z)\right|^p\,d\mu_{u, r} $$
for all $r < 0$. Therefore
$$ \lim_{r\to 0^-} \int_{S_{u,r}}\left|f(z)\right|^pd\mu_{u, r} \le 1,$$
showing that $f\in B_{u,p}(1)$.
\end{proof}

\par Denote by $\E_1$ the family of $u\in \E$ such that $\il \aD \Dl u = 1$ and for such $u$ define 
\begin{align*}
B_{u,p}(R)  & =\{f\in H^p_u(\aD): \|f\|_{H^p_u}\le R\} \text{ and }\\
B_\infty(R) & = \{f\in H^\infty(\aD):|f| \le R\}. 
\end{align*}
Also let $\tilde{\E}_1\sbs \E_1$ consist of those $u\in \E_1$ for which $\al(e^{i\th}) =\il {\bd\aD}P(z, e^{i\th})\,\Dl u(z) <\infty$ for all $\th \in [0,2\pi]$.

\bT \[\bigcap_{u\in \tilde{\E}_1}B_{u,p}(1) = B_\infty(1).\]\eT
\begin{proof}
The inclusion $B_\infty(1) \sbs \bigcap_{u\in \E_1}B_{u,p}(1)$ is clear. For the other way around, let $f\in H^\infty(\aD)\setminus B_\infty(1)$. Since $|f^*(e^{i\th})|^p \in L^1(\lm)$, by the Fatou's theorem 
\[\int_{\bd \aD}P(e^{i\th}, re^{i\hvi})|f^*(e^{i\th})|^p\,d\lm \to |f^*(e^{i\hvi})|^p\, \text{ a.e. }\]
Hence there exists  $A\sbs\{\th\in[0, 2\pi]: f^*(e^{i\th}) \text{ exists}\}$  with  $\lm(A) >0$ such that 
\begin{itemize}
\item $|f^*(e^{i\hvi})|> 1$  and 
\item $\int_{\bd \aD}P(e^{i\th}, re^{i\hvi})|f^*(e^{i\th})|^p\,d\lm \to |f^*(e^{i\hvi})|^p$ 
\end{itemize}
for every $\hvi \in A$. We may suppose that $0\in A$
 
\par Since $u(z) = \int_{\aD} G(z, w)\,\Dl u(w)$, where $G(z,w)$ is the Green's function for the unit disk, and $\frac{\bd }{\bd n}G(z,w)|_{z=e^{i\th}}=P(e^{i\th}, w)$, 
\[\frac{\bd u}{\bd n}(e^{i\th}) = \int_\aD P(e^{i\th}, w)\,\Dl u(w).\]
From section \ref{bv}, for $f\in H^p_u(\aD),\,p>1$,
\[\|f\|^p_{H^p_u} = \int_{\bd \aD}\frac{\bd u}{\bd n}(e^{i\th}) |f^*(e^{i\th})|^p\,d\lm.\]

\par Let $t_k\nea 1$ and $u_k(z) = G(z, t_k)$.  Then
\begin{equation*}\begin{aligned}
\|f\|^p_{H^p_{u_k}} & = \int_{\bd \aD}P(e^{i\th}, t_k)|f^*(e^{i\th})|^p\,d\lm\\
&\longrightarrow |f^*(1)|^p\, \text{ as } k\to \infty \text{ because } 0\in A. 
\end{aligned}\end{equation*}
 Hence $f\not\in \bigcap_{u\in\tilde{\E}_1} B_{u,p}(1).$ The theorem follows.
\end{proof}

\par Recall from Proposition \ref{abscont} that we have $\mu_u = \al \lm$ where $\al\in L^1(\lm)$ and $\al \ge c >0$ for some constant $c$. Moreover, $\al$ is lower semicontinuous. Hence, there exists an increasing sequence of positive smooth functions $\al_n$ converging to $\al$ pointwise. Define
\[\tilde{\al}(z)  = \il {\aT} \frac{e^{i\th} + z}{e^{i\th} - z}\log\al(e^{i\th})\,d\lm(\th)\]
\[\tilde{\al}_n(z) = \il {\aT}\frac{e^{i\th}+z}{e^{i\th} - z}\log\al_n(e^{i\th})\,d\lm(\th).\]
Clearly $\tilde{\al},\,\tilde{\al_j}\in \rk{O(\aD)}$, so the functions $A(z) = e^{\tilde{\al}(z)}\text{ and } A_n(z) = e^{\tilde{\al}_n(z)}\in \rk{O(\aD)}$. Moreover, The functions $\tilde\al_n$ and $A_n$ extend smoothly to the boundary, $|A^*(e^{i\th})|= \al(e^{i\th})$ and $|A_n^*(e^{\th})|= \al_n(e^{i\th})$.
\bT\label{isom} The space $H^p_u(\aD)$ is isometrically isomorphic to $H^p(\aD)$. \eT
\begin{proof}
First, we show that if $f\in H^p_u(\aD)$ then $ A^{1/p}f\in H^p(\aD)$. Clearly $A_n^{1/p}f \in H^p(\aD)$. Then by formula (9) in \cite[IX.4]{Goluzin},
\[\int_0^{2\pi}|A_n(re^{i\th})||f(re^{i\th})|^p\,d\lm(\th) \le \int_0^{2\pi} |A_n^*(e^{i\th})||f^*(e^{i\th})|^p\,d\lm(\th).\]
 Since $A_n^{1/p}f$ converges to $A^{1/p}f$ uniformly on compact subsets of $\aD$, for $0< r < 1$,
\begin{align*}
\int_0^{2\pi}|A(re^{i\th})||f(re^{i\th})|^p\,d\lm & =\lim_{n\to\infty} \int_0^{2\pi}|A_n(re^{i\th})||f(re^{i\th})|^p\,d\lm(\th)\\
& \le \lim_{n\to\infty} \int_0^{2\pi} |A_n^*(e^{i\th})||f^*(e^{i\th})|^p\,d\lm(\th)\\
& = \|f\|^p_{H^p_u}.
\end{align*}
The last equality above follows from the monotone convergence theorem. Thus $A^{1/p}f \in H^p(\aD)$.
\par
Now, define an operator
\begin{align*}
\Phi:H^p_u(\aD) & \to H^p(\aD)\\
   f & \mapsto A^{1/p}f.
\end{align*}
Clearly $\Phi$ is linear.  Since
\[\int_0^{2\pi} |A^*(e^{i\th})||f^*(e^{i\th})|^p\,d\lm = \int_0^{2\pi} |f^*(e^{i\th})|^p\al(e^{i\th})\,d\lm=\il {\aT}|f^*|^p\,d\mu_u,\]
we have $\|A^{1/p}f\|_{H^p} =\|f\|_{H^p_u}$. So $\Phi$ is an isometry.
\par Let $f\in H^p(\aD)$. Since $|A(z)| \ge c > 0$,  $A^{-1/p}f \in H^p(\aD)$. It follows from the identity
\[\il {\aT} |A^*|^{-1}|f^*|^p\,d\mu_u  = \il {\aT}|f^*|^p\,d\lm\]
together with Theorem \ref{bvalue} that $A^{-1/p}f \in H^p_u(\aD)$.
Thus $\Phi$ is a surjective linear isometry. We are done.
\end{proof}

\end{document}